\documentclass[conference,compsoc]{IEEEtran}
\ifCLASSOPTIONcompsoc
  \usepackage[nocompress]{cite}
\else
  \usepackage{cite}
\fi
\ifCLASSINFOpdf
  \usepackage[pdftex]{graphicx}
   \graphicspath{{../pdf/}{../jpeg/}}
   \DeclareGraphicsExtensions{.pdf,.jpeg,.png}
\else
   \usepackage[dvips]{graphicx}
   \graphicspath{{eps/}}
   \DeclareGraphicsExtensions{.eps}
\fi
\usepackage{amsmath}
\ifCLASSOPTIONcompsoc
  \usepackage[caption=false,font=footnotesize,labelfont=sf,textfont=sf]{subfig}
\else
  \usepackage[caption=false,font=footnotesize]{subfig}
\fi
\begin{document}
%
\title{Local Linear Constraint based Optimization Model for Dual Spectral CT}

\author{\IEEEauthorblockN{Qian Wang}
\IEEEauthorblockA{Department of Electrical and Computer Engineering\\
University of Massachusetts Lowell\\
Lowell, Massachusetts 01854\\
Email: Qian\_Wang@student.uml.edu}}

\maketitle

\begin{abstract}
Dual spectral computed tomography (DSCT) can achieve energy- and material-selective images, and has a superior distinguishability of some materials than conventional single spectral computed tomography (SSCT). However, the decomposition process is illposed, which is sensitive with noise, thus the quality of decomposed images are usually degraded, especially the signal-to-noise ratio (SNR) is much lower than single spectra based directly reconstructions. In this work, we first establish a local linear relationship between dual spectra based decomposed results and single spectra based directly reconstructed images. Then, based on this constraint, we propose an optimization model for DSCT and develop a guided image filtering based iterative solution method. Both simulated and real experiments are provided to validate the effectiveness of the proposed approach.
\end{abstract}


%
\IEEEpeerreviewmaketitle

\section{Introduction}
In X-ray dual spectral computed tomography (DSCT), a specimen is scanned with two different X-ray energy spectra, then the collected polychromatic projections from this procedure are utilized to perform energy- and material-selective reconstruction \cite{alvarez1976energy,alvarez1979comparison,kalender1986evaluation,vetter1986evaluation,kalender1988algorithm,chuang1988comparison}. Compared with conventional single spectral computed tomography (SSCT), DSCT has many advantages for contrast enhancement, artifact reduction, and material composition assessment. Thus, DSCT has wide potential applications in medical and industrial domains, such as bone mineral density and liver iron concentrations measurements, beam-hardening correction and contrast enhancement of soft tissue, positron emission tomography (PET) attenuation correction, rock core characterization for petrochemical industry, and so forth \cite{coleman1985beam,fessler2002maximum,zhang2006exact,kinahan2006dual,ying2006dual,zhang2008practical,noh2009statistical,johnson2011dual}.

Existing methods to perform the decomposition of DSCT can be classified into three groups: image based methods, projection based methods, and iterative methods. Image based methods treat the projection data sets as being independent until they are reconstructed. Then, images of each spectra are combined linearly to obtain two decomposed images \cite{maass2009image}. Such methods fail to describe the real nonlinearity relationship between decomposed results and polychromatic projections, thus beam hardening artifacts cannot be removed successfully \cite{brooks1976beam,coleman1985beam}. Projection based methods correctly treat the available information by passing the projection data through a high order decomposition function, followed by image reconstruction \cite{flohr2006first,stenner2007empirical}. Generally, they can obtain better decomposition performances than image based ones, but the combination of polychromatic projections before reconstruction requires a satisfication of geometrical consistency. Several iterative methods are proposed based on statistical model and nonlinear optimization \cite{elbakri2002statistical,xu2009implementation,maass2009dual,niu2014iterative}. By introducing prior knowledge or establishing nonlinear model, these methods improve the quality of decomposed results effectively, but they usually have a drawback of slow convergence and huge computational cost. An extended algebraic reconstruction technique (E-ART) for DSCT was proposed by Zhao et al. \cite{zhao2015extended}. It describes the DSCT reconstruction as a nonlinear system problem, and then extends the classic ART method to solve the model iteratively. This method can produce high quality decompositions from polychromatic projections, but it also has the drawback of slow convergence. Hu et al. developed the E-ART method into an extended simultaneous algebraic reconstruction version, i.e. E-SART \cite{hu2016extended}. This method is based on the matrix inversion and has a high degree of parallelism, thus the convergence rate is improved dramatically. But the illposedness of the decomposition process render it noticeably sensitive with noise, thus the achieved quality is degraded with reduced signal-to-noise ratio (SNR).

Comparing decomposed results of DSCT with reconstructed images of SSCT, we find that although single spectra based reconstructions have weaker capability of destinguishing some materials, the achieved quality like SNR is dramatically higher. Moreover, there is an interesting relationship between them, i.e., decomposed results of DSCT can be considered as modifications of reconstructed images of SSCT by removing some components and adjusting gray values. Further, this structure-based feature can be described as a local linear relationship mathematically. By combing it as a constraint into an optimization model for DSCT, the reconstructed image of SSCT works as a guided image and contributes to improving the smoothness of decomposed results of DSCT effectively. Thus, the systematic noise are well surpressed and the quality of decomposed results are improved considerably.

In this work, a local linear model is first established to describe the relationship between decomposed results of DSCT and reconstructed images of SSCT. Then by employing it as a constraint, we propose an optimization model for DSCT. Based on the guided image filtering \cite{he2010guided,he2013guided}, a correlative iterative solution method is developed, which produces decomposed images with high SNR and maintains a fast convergence meanwhile.

The remainder of this paper is organized as follows. In section \ref{sec:methods}, the mathematical model of DSCT is presented, then the E-SART method and the guided image filtering technique are reviewed briefly. In section \ref{sec:work}, we first propose a local linear constraint based optimization model DSCT. Then we develop a guided image filtering based iterative solution method. In section \ref{sec:experiments}, both simulated and real experiments are provided to verify the effectiveness of the proposed methods. Section \ref{sec:conclusion} contains our conclusions.

\section{Methods}
\label{sec:methods}
\subsection{Mathematical Model of DSCT}
By considering the X-ray is polychromatic and assuming collected raw data are consistent geometrically, we can describe the physical process of DSCT as follows,
\small
\begin{equation}
\label{eq:poly}
	P_{k,l}=-\ln\int_{E}S_{k}(E)\exp{\big(-\mathcal{P}_{l}(\mu(E,\mathbf{x}))\big)} \,\mathrm{d}E, l\in\mathcal{L}, k=1,2,
\end{equation}
\normalsize
where $\mu(E,\mathbf{x})$ is the linear attenuation coefficient of a specimen at a spatial position $\mathbf{x}$ and energy $E$; $\mathcal{P}_{l}(\cdot)$ represents the ray transform which is an integral transform along a ray path $l$; $S_{k}(E)$ is the $k$-th normalized emission spectrum; $P_{k,l}$ indicates the acquired information, frequently called projection data.

In DECT, the linear attenuation coefficient $\mu(E,\mathbf{x})$ is usually considered splittable with variable $E$ and $\mathbf{x}$, i.e.,
\begin{equation}
\label{eq:decom}
	\mu(E,\mathbf{x}))=\sum_{i=1}^{2}\psi_{i}(E)f_{i}(\mathbf{x}),
\end{equation}
where $\psi_{i}(E)$ is a function of energy; $f_{i}(\mathbf{x})$ is a function of spatial position. There are commonly two physical explanations of eq. \eqref{eq:decom}: basis material based decomposition and effect based decomposition. For the former, $\psi_{i}(E)$ is the mass attenuation coefficient for material $i$, and $f_{i}(\mathbf{x})$ represents the correlative density distribution. For the latter, $\psi_{1}(E)=E^{-3}$ and $\psi_{2}(E)=KN(E)$ (Klein-Nishina function) corresponds to the photoelectric effect and Compton scattering respectively, and $f_{i}(\mathbf{x})$ represents the effect distribution correspondingly. The aim of DSCT is to reconstruct images of distrubution functions $f_{i}(\mathbf{x}), i=1,2$.

\subsection{E-SART Method}
\label{subsec:ESART}
By substituting eq. \eqref{eq:decom} into eq. \eqref{eq:poly} and discretizing the correlative result according to energy bin, we get
\begin{equation}
\label{eq:disc}
	P_{k,l}=-\ln\sum_{j=1}^{J_{k}}S_{k,j}\exp{\big(-\sum_{i=1}^{2}\psi_{i,j}\mathcal{P}_{l}(\mathbf{f}_{i})\big)} \Delta E,
\end{equation}
where $J_{k}$ is the enery bin number of spectra $k$; $\Delta E$ represents the bin length; $S_{k,j}$ and $\psi_{i,j}$ are the samplings of $S_{k}(E)$ and $\psi_{i}(E)$ within bin $j$; $\mathbf{f}_{i}$ is a one demensional column vector representing the discretized distribution function. Then, the 1st order Taylor expansion of eq. \eqref{eq:disc} at point $(\mathbf{f}_{1}(n);\mathbf{f}_{2}(n))$ is
\begin{equation}
\label{eq:taylor}
P_{k,l}\approx P_{k,l}(n)+\begin{pmatrix}\frac{\Psi_{k,l}^{1}(n)}{Q_{k,l}(n)},\frac{\Psi_{k,l}^{2}(n)}{Q_{k,l}(n)}\end{pmatrix}\begin{pmatrix}\mathcal{P}(\mathbf{f}_{1}-\mathbf{f}_{1}(n))\\ \mathcal{P}(\mathbf{f}_{2}-\mathbf{f}_{2}(n))\end{pmatrix},
\end{equation}
where
\setlength{\arraycolsep}{0.0em}
\begin{eqnarray*}
P_{k,l}(n)&=&-\ln\sum_{j=1}^{J_{k}}S_{k,j}\exp{\big(-\sum_{i=1}^{2}\psi_{i,j}\mathcal{P}_{l}(\mathbf{f}_{i}(n))\big)} \Delta E,\\
Q_{k,l}(n)&=&\sum_{j=1}^{J_{k}}S_{k,j}\exp{\big(-\sum_{i=1}^{2}\psi_{i,j}\mathcal{P}_{l}(\mathbf{f}_{i})\big)} \Delta E,\\
\Psi_{k,l}^{1}(n)&=&\sum_{j=1}^{J_{k}}\psi_{1,j}S_{k,j}\exp{\big(-\sum_{i=1}^{2}\psi_{i,j}\mathcal{P}_{l}(\mathbf{f}_{i})\big)} \Delta E,\\
\Psi_{k,l}^{2}(n)&=&\sum_{j=1}^{J_{k}}\psi_{2,j}S_{k,j}\exp{\big(-\sum_{i=1}^{2}\psi_{i,j}\mathcal{P}_{l}(\mathbf{f}_{i})\big)} \Delta E.
\end{eqnarray*}
\setlength{\arraycolsep}{5pt}
Along each ray path, two projection data are acquired based on different X-ray spectra, i.e., low-energy and high-energy. By solving the system of linear equation \eqref{eq:taylor}, $k=1,2$, we can get the projection of distribution fuction in a iteration form,
\scriptsize
\begin{equation}
\begin{pmatrix}\mathcal{P}(\mathbf{f}_{1}(n+1))\\ \mathcal{P}(\mathbf{f}_{2}(n+1))\end{pmatrix}=\begin{pmatrix}\mathcal{P}(\mathbf{f}_{1}(n))\\ \mathcal{P}(\mathbf{f}_{2}(n))\end{pmatrix}+\frac{C_{l}(n)}{\det(M_{l}(n))}\begin{pmatrix}P_{1,l}-P_{1,l}(n)\\ P_{2,l}- P_{2,l}(n)\end{pmatrix},
\end{equation}
\normalsize
where
\begin{equation*}
	M_{l}(n)=\begin{pmatrix}\frac{\Psi_{1,l}^{1}(n)}{Q_{1,l}(n)},\frac{\Psi_{1,l}^{2}(n)}{Q_{1,l}(n)}\\\frac{\Psi_{2,l}^{1}(n)}{Q_{2,l}(n)},\frac{\Psi_{2,l}^{2}(n)}{Q_{2,l}(n)}\end{pmatrix},
	C_{l}(n)=\begin{pmatrix}\frac{\Psi_{2,l}^{2}(n)}{Q_{2,l}(n)},-\frac{\Psi_{1,l}^{2}(n)}{Q_{1,l}(n)}\\-\frac{\Psi_{2,l}^{1}(n)}{Q_{2,l}(n)},\frac{\Psi_{1,l}^{1}(n)}{Q_{1,l}(n)}\end{pmatrix}.
\end{equation*}
Then by using the conventional Simultaneous Algebraic Reconstruction Technique (SART), distribution fuction $\mathbf{f}_{1}$ and $\mathbf{f}_{2}$ are updated iteratively.

Comparing with E-ART, E-SART improves the convergence rate dramatically. However, the illposedness of the inverse problem render this matrix inversion based decomposition process sensitive with inevitable systematic noise, i.e., obtained results suffer from low SNR. Thus, some prior knowledge or constraints are needed to improve the robust against noise.

\subsection{Guided Image Filtering}
\label{subsec:guided}
Guided filter is an edge-preserving filter with a great variety of applications \cite{levin2004colorization,levin2008closed,farbman2008edge,he2011single}, of which the key assumption is a local linear model between the guidance image $I$ and the filtering output $y$, i.e.,
\begin{equation}
\label{eq:linear}
	y_{i}=a_{k}I_{i}+b_{k}, \forall i\in\omega_{k},
\end{equation}
where $\omega_{k}$ is a window centered at the pixel $k$ with radius $r$; $(a_{k}, b_{k})$ are some linear coefficients constant in $\omega_{k}$. Modeling the output $y$ as the input $x$ removing some unwanted noise or textures $t$:
\begin{equation*}
	y_{i}=x_{i}-t_{i}.
\end{equation*}
Thus, by minimizing the difference between $y$ and $x$ within a window $\omega_{k}$ while maintaining the linear model \eqref{eq:linear}, the correlative optimization model is established as follows,
\begin{equation*}
	\min\sum_{i\in\omega_{k}}\big((a_{k}I_{i}+b_{k}-x_{i})^{2}+\epsilon a_{k}^{2}\big),
\end{equation*}
where $\epsilon$ is a regularization parameter penalizing large $a_{k}$. The solution is given by
\setlength{\arraycolsep}{0.0em}
\begin{eqnarray*}
	a_{k}&=&\frac{\frac{1}{|\omega|}\sum_{i\in \omega_{k}}I_{i}x_{i}-\nu_{k}\bar{x}_{k}}{\sigma_{k}^{2}+\epsilon},\\
	b_{k}&=&\bar{x}_{k}-a_{k}\nu_{k},
\end{eqnarray*}
\setlength{\arraycolsep}{5pt}
where $\nu_{k}$ and $\sigma_{k}^{2}$ are the mean and variance of $I$ in $\omega_{k}$; $|\omega|$ is the number of pixels in $\omega_{k}$; $\bar{x}_{k}$ is the mean of $x$ in $\omega_{k}$. Then, the filtering output $y$ can be computed according to eq. \eqref{eq:linear}. Because a pixel $i$ is involved in all the covered windows, by averagig all the possibles output values, we get
\begin{equation*}
		y_{i}=\bar{a}_{i}I_{i}+\bar{b}_{i}.
\end{equation*}
Here $\bar{a}_{i}$ and $\bar{b}_{i}$ are the average coefficients of all windows overlapping $i$. By using the guided image filtering, input $x$ is refined by the guided image $I$ according to the local linear relationship between them.

\section{Local Linear Constraint based Optimization Model and Iterative Solution Method}
\label{sec:work}
\subsection{Local Linear Constraint based Optimization Model}
\begin{figure*}[htbp]
	\centering
	\includegraphics[width=7in]{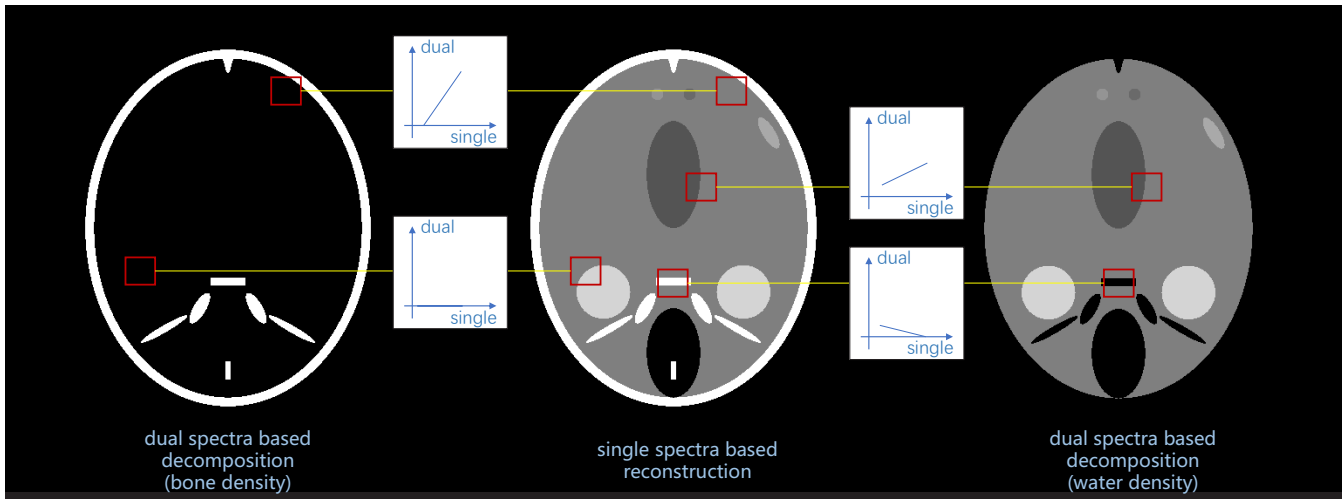}
	\caption{An illustration of the local linear relationship: single spectra based directly reconstructed image is shown in the middle column; dual spectra based decomposed results are shown in the left (bone density) and the right (water density) columns.}
	\label{fig:illustration}
\end{figure*}

Although directly reconstructed images under single spectra scan have weak capability of distinguishing some materials, the quality is improved remarkbaly than dual spectra based decomposition results, especially when the noise level is relatively high. Moreover, there are some structure-based relationships between them.  An intuitionistic character is that decomposed results of DSCT can be viewed as modifications of reconstructed images of SSCT by removing some components and adjusting gray values. When analyzing this feature in detail, as is illustrated in Fig. \ref{fig:illustration}, we can find that a linear relationship is usually held in small patches, of which the discreted version reads,
\begin{equation*}
	\mathbf{f}^{(j)}=a^{(k)}\mathbf{g}^{(j)}+b^{(k)}, \forall j\in\omega^{(k)},
\end{equation*}
where $\mathbf{f}$ represnets a decomposed result of DSCT; $\mathbf{g}$ represents a reconstructed image of SSCT; $j$ is a pixel index; $\omega^{(k)}$ is a window centered at the pixel $k$ with radius $r$; $(a^{(k)}, b^{(k)})$ are some linear coefficients constant in $\omega^{(k)}$.

Based on this constraint, we proposed an optimization model for DSCT as follows,
\setlength{\arraycolsep}{0.0em}
\scriptsize
\begin{eqnarray}
\label{eq:opt}
\min_{(\mathbf{f}_{1},\mathbf{f}_{2})}\bigg\{\Big\|\begin{pmatrix}\mathcal{P}(\mathbf{f}_{1})-\mathcal{P}(\mathbf{f}_{1}(n))\\ \mathcal{P}(\mathbf{f}_{2})-\mathcal{P}(\mathbf{f}_{2}(n))\end{pmatrix}-\frac{C_{l}(n)}{\det(M_{l}(n))}\begin{pmatrix}P_{1,l}-P_{1,l}(n)\\ P_{2,l}- P_{2,l}(n)\end{pmatrix}\Big\|_{L_{2}}^{2}\nonumber\\
+\sum_{i=1}^{2}\sum\limits_{k}\sum\limits_{j\in\omega_{i}^{(k)}}\frac{\xi_{i}}{|\omega_{i}|}\Big(\big(a_{i}^{(k)}\mathbf{g}_{\mathbf{f}_{i}}^{(j)}+b_{i}^{(k)}-\mathbf{f}_{i}^{(j)}(n)\big)^{2}+\epsilon_{i}(a_{i}^{(k)})^{2}\Big)\bigg\},
\end{eqnarray}
\normalsize
\setlength{\arraycolsep}{5pt}
where $\mathbf{g}_{\mathbf{f}_{i}}$ is a selected result of SSCT corresponding to $\mathbf{f}_{i}$; $\xi_{i}$ and $\epsilon_{i}$ are regularization parameters.

In model \eqref{eq:opt}, for each searched-for decomposition result $\mathbf{f}_{i}, i=1,2$, we employ a correlative local linear constraint. Thus, the smoothness knowledge of reconstructed image of SSCT is effectively introduced into the decomposition problem. By  weakening the illposedness, the noise is surpressed noticeably, thus the quality of decomposed result, like SNR, is improved dramatically.

\subsection{Iterative Solution Method}
When solving the proposed optimization model, considering the data term is measured in projection domain and the regularity terms are measured in image domain, we split model \eqref{eq:opt} into two sub-optimization problems and develop an iterative scheme as follows,
\scriptsize
 \begin{subequations}
\begin{align}
\label{eq:opt2_pfstar}
&\big(\mathcal{P}(\mathbf{f}_{1}(*)),\mathcal{P}(\mathbf{f}_{2}(*))\big)=\nonumber\\
&\min_{\big(\mathcal{P}(\mathbf{f}_{1}),\mathcal{P}(\mathbf{f}_{2})\big)}\Big\|\begin{pmatrix}\mathcal{P}(\mathbf{f}_{1})-\mathcal{P}(\mathbf{f}_{1}(n))\\ \mathcal{P}(\mathbf{f}_{2})-\mathcal{P}(\mathbf{f}_{2}(n))\end{pmatrix}-\frac{C_{l}(n)}{\det(M_{l}(n))}\begin{pmatrix}P_{1,l}-P_{1,l}(n)\\ P_{2,l}- P_{2,l}(n)\end{pmatrix}\Big\|_{L_{2}}^{2},\\
\label{eq:opt2_fstar}
&\big(\mathbf{f}_{1}(*),\mathbf{f}_{2}(*)\big)=\mathcal{P}^{-1}\big(\mathcal{P}(\mathbf{f}_{1}(*)),\mathcal{P}(\mathbf{f}_{2}(*))\big),\\
\label{eq:opt2_ab}
&\big(a_{i}(n+1)),b_{i}(n+1))\big)=\nonumber\\
&\min_{(a_{i},b_{i})}\sum\limits_{k}\sum\limits_{j\in\omega_{i}^{(k)}}\frac{\xi_{i}}{|\omega_{i}|}\Big(\big(a_{i}^{(k)}\mathbf{g}_{\mathbf{f}_{i}}^{(j)}+b_{i}^{(k)}-\mathbf{f}_{i}^{(j)}(*)\big)^{2}+\epsilon_{i}(a_{i}^{(k)})^{2}\Big),\\
\label{eq:opt2_abar}
&\bar{a}_{i}^{(k)}=\frac{1}{|\omega_{i}|}\sum\limits_{j\in\omega_{i}^{(k)}}a_{i}^{(j)},\\
\label{eq:opt2_bbar}
&\bar{b}_{i}^{(k)}=\frac{1}{|\omega_{i}|}\sum\limits_{j\in\omega_{i}^{(k)}}b_{i}^{(j)},\\
\label{eq:opt2_f}
&\mathbf{f}_{i}(n+1)=\bar{a}_{i}\mathbf{g}_{\mathbf{f}_{i}}+\bar{b}_{i},
\end{align}
\end{subequations}
\normalsize
where Eqs. \eqref{eq:opt2_ab}-\eqref{eq:opt2_f} are implemented respectively for $i=1,2$. We use the E-SART method reviewed in \ref{subsec:ESART} to solve Eqs. \eqref{eq:opt2_pfstar} and \eqref{eq:opt2_fstar}. When solving Eqs. \eqref{eq:opt2_ab}-\eqref{eq:opt2_f}, the guided image filtering reviewed in \ref{subsec:guided} is employed.

For one thing, the proposed method suppresses the noise effectively, thus the quality of decomposed results is improved dramatically. For another, it maintains the merit of high convergence rate of the E-SART meanwhile. By considering the relatively weaker nonlinearity of high-energy spectra based data than low-energy spectra based ones, in this study we just use the high-energy spectra based reconstructed image as the guided image for both decomposition results $f_{i}, i=1,2$. More analyses and comparisons are still needed to choose a more satisfactory guided image, which will be fully studied in our continuous work. 

\section{Experiments}
\label{sec:experiments}
In order to verify the effectiveness of the proposed method, experiments with both simulated and real data are performed. In the simulated experiments, both noise-free case and noisy case are tested. These experiments are restricted to
fan beam CT for simplicity. The generalization to cone beam case is straightforward. As a comparison, we implement the E-SART method as well, because it is a relatively new method with fast convergence rate. We employ all the methods for the case of basis material based decomposition, which can be extended to the case of effect based case easily. And the iteration number is fixed to 30.

\subsection{Simulated Experiments}
\begin{figure}[htbp]
	\centering
	\includegraphics[width=3.6in]{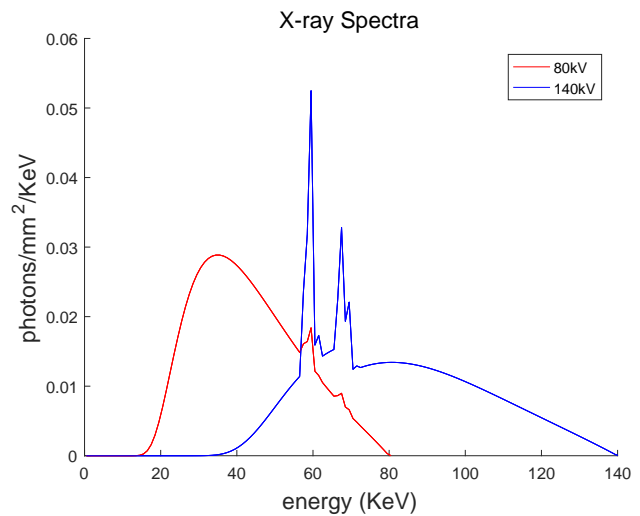}
	\caption{The X-ray spectra used in the simulated experiments.}
	\label{fig:spectra}
\end{figure}

The phantom used in all the simulated experiments is the 2D FORBILD head phantom without ears shown in Fig. \ref{fig:illustration} \cite{phantom}. Water and bone are chosen as the two basis materials, of which the mass attenuation coefficients are retrieved from the National Institute of Standard Technology (NIST) tables of X-ray mass attenuation coefficient \cite{hubbell1995tables}. The polychromatic spectra of a GE Maxiray 125 X-ray tube are simulated by using an open source X-ray spectra simulator, SpectrumGUI \cite{SpectrumGUI}. Two tube voltages, 80 kV and 140 kV, are chosen, where the latter is filtered with 1mm copper. The correlative spectra are shown in Fig. \ref{fig:spectra}. The energy of photons emitted from the source is 8 MeV. The detector consists of 512 channels with length 0.03 cm. The source-object distance (SOD) is 100 cm and the source-detector distance (SDD) is 120 cm. With this configuration, images are reconstructed to a 512 $\times$ 512 digitization with a pixel size of 0.0249 cm.

\begin{figure*}[htbp]
	\centering
	\includegraphics[width=7in]{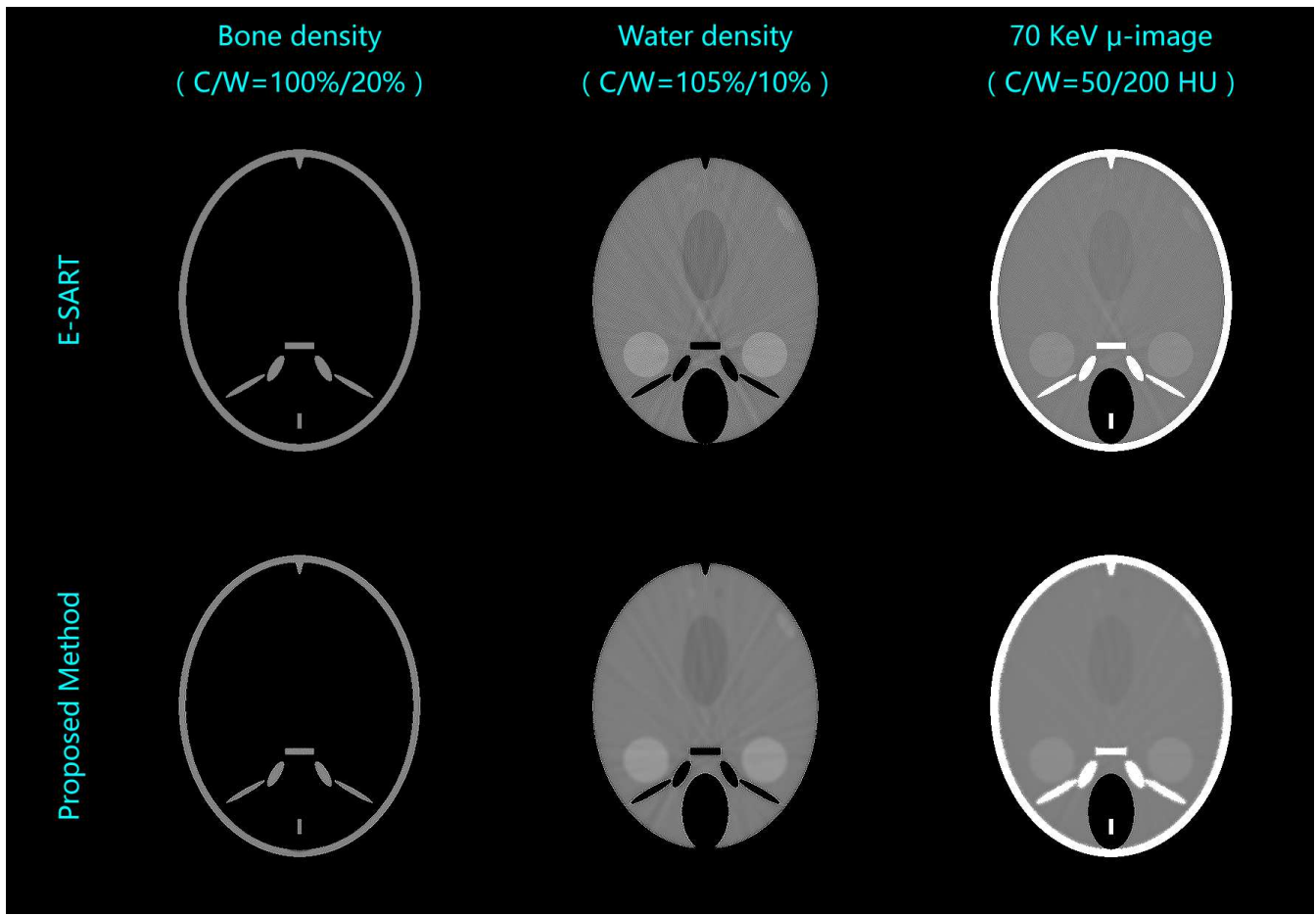}
	\caption{Simulated experiments of the noise-free case. The results of the E-SART method and the proposed method are shown in the first row and the second row respectively. The bone and water based density images are presented in the first two columns. The composite images at 70 KeV are shown in the last column.}
	\label{fig:simulated_I}
\end{figure*}

\begin{figure*}[htbp]
	\centering
	\includegraphics[width=7in]{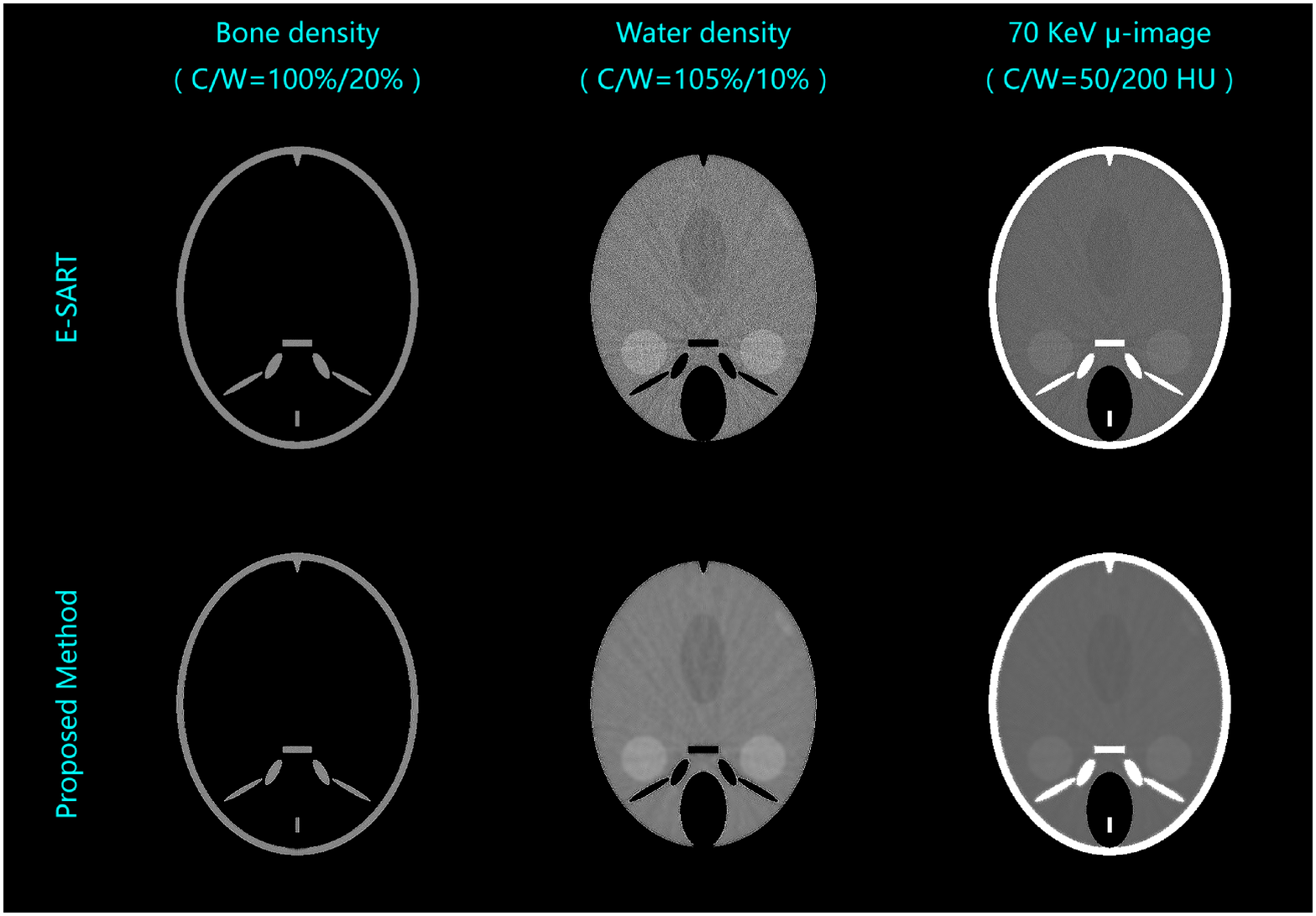}
	\caption{Simulated experiments of the noisy case. The results of the E-SART method and the proposed method are shown in the first row and the second row respectively. The bone and water based density images are presented in the first two columns. The composite images at 70 KeV are shown in the last column.}
	\label{fig:simulated_N}
\end{figure*}

Under this setting, we test both the noise-free case and the noisy case, seen Fig. \ref{fig:simulated_I} and Fig. \ref{fig:simulated_N} respectively. All the experimenal resutls demonstrate the conspicuous noise tolerance of the proposed method.

\subsection{Real Experiments}
In the real experiment, the measured specimen is a bone submerged in water. An X-ray source (YXLON Y.TU450 D09 tube) is operated at the tube voltage of 80 kV and 140 kV for low- and high-energy spectra scan respectively, and the tube current is 5 mA. The employed flat-pannel detector (YXLON Y.LDA detector) has $1920\times1920$ channels with mesh size of 0.0127 cm. The SOD is 23.15 cm and the SDD is 69.67 cm. By using a collimator, the data from the central slice are obtained to validate the proposed method. The iteration number is 10.

\begin{figure*}[htbp]
\centering
\includegraphics[width=7in]{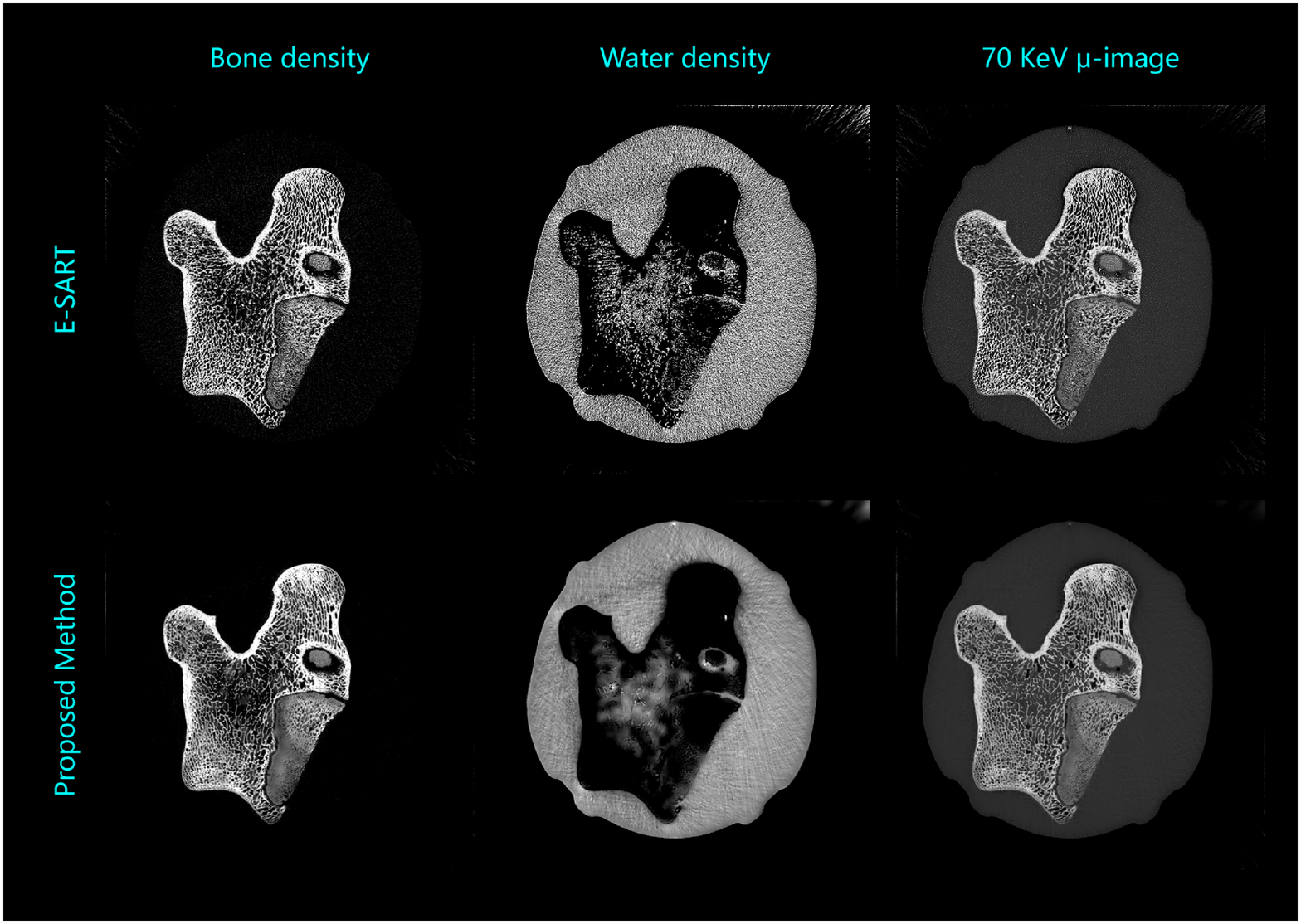}
\caption{Real experiments by using the E-SART method (the first row) and the proposed method (the second row). The bone and water based density images are presented in the first two columns respectively. The composite images at 70 KeV are shown in the last column.}
\label{fig:real}
\end{figure*}

The decomposed results by using the E-SART method and the proposed method are shown in Fig. \ref{fig:real}. It is noticeable that the proposed method can surpress the noice effectively and improve the smoothness dramatically.

\section{Conclusion}
\label{sec:conclusion}
In this work, we first establish a local linear constraint to describe the structure relationship between dual spectra based decomposed results and single spectra based directly reconstructed images. Then the correlative optimization model and the iterative solution method are proposed respectively. By employing the guided image filtering, the smoothness knowledge of the single spectra based reconstructed image is effectively introduced into the decomposition process. This method reduces the illposedness of the DSCT, thus the noise in the decomposition process is surpressed successfully. Both simulated and real experiments demonstrate the effectiveness of the proposed method. Further studies are needed to select a satisfactory guided image.






%

\bibliographystyle{IEEEtran}
\bibliography{Reference}

\end{document}